\documentclass[11pt,twoside,leqno]{aomamlt2e} 
  \pageno{1097}
\received{October 13, 2004}

\numberwithin{equation}{section}

\theoremstyle{definition}

\theoremstyle{remark}

\begin{document}
\currannalsline{164}{2006} 

 \title{Erratum to ``Generalizations of\\ the Poincar\'e-Birkhoff theorem''}

 \acknowledgements{Supported in part by NSF grant DMS-0099640.}
\author{John Franks}

  \institution{Northwestern University, Evanston, IL
\\
 \email{john@math.northwestern.edu}}

 \shorttitle{Generalizations of the Poincar\'e-Birkhoff theorem}

\def\ti{\tilde}
\def\sinfty{S_{\infty}}
\def\sl3z{SL(3, \mathbb Z)}
\def\ti{\tilde}
\def\rtwo{\mathbb R^2}
\def\rone{\mathbb R}
\def\tag{translation arc geodesic}
\def\nc{near cycle}
\def\calo{{\cal O}}
\def\calb{{\cal B}}
\def\D{{\cal D}}
\def\G{{\cal G}}
\def\H{{\cal H}}
\def\A{{\mathbb A}}
\def\Q{{\mathbb Q}}
\def\RH{RH(W,\partial_+ W)}
\def\nr{{\bf NEED REFERENCE\ }}

Theorem (2.1) of \cite{Fr} is not stated correctly.  The correct statement
(and the result which is actually proved) is the following.

\demo{\scshape Theorem}
{\it Suppose $f:\A\to \A$ is an orientation preserving 
homeomorphism of the open annulus which is homotopic to the identity{\rm ,} and 
satisfies the following conditions\/{\rm :}\/
\begin{itemize}
\ritem{1.} Every point of $\A$ is non-wandering.

\ritem{2.} $f$ has at most finitely many fixed points.

\ritem{3.} There is a lift of $f$ to its universal covering space{\rm ,} $\ti
f:\ti \A\to \ti \A${\rm ,} which possesses both a positively returning disk{\rm ,}
which is a lift of a disk in $\A${\rm ,} and a negatively returning disk{\rm ,}
which is a lift of a disk in $\A$.
\end{itemize}
Then $f$ has a fixed point of positive index. }
\Enddemo

Other results of \cite{Fr} which rely on Theorem (2.1) are correct as
stated and follow from the result above.  

A disk $U$ in the covering space $\ti \A$ is called {\it positively
returning} for $\ti f$ if $\ti f(U)\cap U=\emptyset$ and $\ti
f^n(U)\cap T^k(U)\ne\emptyset$ for some $n,k>0$, where $T$ is the
generator of the infinite cyclic group of covering translations.
Negatively returning disks are defined similarly.  The difference in
the statement of the theorem above and the statement of
Theorem (2.1) of \cite{Fr} is
the additional requirement in item 3 that the positively and
negatively returning disks in the covering space are lifts of disks in
the annulus.  Equivalently the positively and negatively returning
disks in the covering space must be disjoint from their image under
$T^k,\ k\ne 0$ (standard arguments show it suffices
that $T(U) \cap U = \emptyset$). This additional hypothesis is
satisfied in other results from \cite{Fr} which use Theorem (2.1).

I do not know whether Theorem (2.1) of \cite{Fr},  as originally stated,
is true or not.


\references{100}

\bibitem[1]{Fr}
\name{J.~Franks,}
\newblock Generalizations of the Poincar\'e-Birkhoff theorem,
\newblock {\it Ann.\ of Math\/}.\ {\bf 128} (1988), 139--151.
\Endrefs
\end{document}